\documentclass[11pt]{amsart}
\usepackage[latin1]{inputenc}
\usepackage{amscd}
\usepackage{latexsym}
\usepackage{pstcol,pst-plot,pst-3d}
\usepackage{mathptmx}
\usepackage{multicol}
\def\NZQ{\Bbb}               
\def\NN{{\NZQ N}}

\def\AA{{\NZQ A}}

\def\frk{\frak}               

\def\mm{{\frk m}}

\def\opn#1#2{\def#1{\operatorname{#2}}} 
\opn\chara{char}
\opn\length{\ell}
\opn\pd{pd}
\opn\rk{rk}
\opn\projdim{proj\,dim}
\opn\injdim{inj\,dim}
\opn\rank{rank}
\opn\depth{depth}
\opn\grade{grade}
\opn\height{height}
\opn\embdim{emb\,dim}
\opn\codim{codim}

\opn\Tr{Tr}
\opn\bigrank{big\,rank}
\opn\superheight{superheight}\opn\lcm{lcm}
\opn\trdeg{tr\,deg}%
\opn\reg{reg}
\opn\lreg{lreg}
\opn\skel{skel}
\opn\com{com}
\opn\div{div}
\opn\Div{Div}
\opn\cl{cl}
\opn\Cl{Cl}
\opn\Spec{Spec}
\opn\Supp{Supp}
\opn\supp{supp}
\opn\Sing{Sing}
\opn\Ass{Ass}
\opn\Ann{Ann}
\opn\Rad{Rad}
\opn\Soc{Soc}
\opn\Ker{Ker}
\opn\Coker{Coker}
\opn\Im{Im}
\opn\Hom{Hom}
\opn\Tor{Tor}
\opn\Ext{Ext}
\opn\End{End}
\opn\Aut{Aut}
\opn\id{id}

\opn\nat{nat}
\opn\pff{pf}
\opn\Pf{Pf}
\opn\GL{GL}
\opn\SL{SL}
\opn\mod{mod}
\opn\ord{ord}
\opn\aff{aff}
\opn\con{conv}
\opn\relint{relint}
\opn\st{st}
\opn\lk{lk}
\opn\cn{cn}
\opn\core{core}
\opn\vol{vol}
\opn\link{link}
\opn\star{star}
\opn\gr{gr}

\def\pot#1#2{#1[\kern-0.28ex[#2]\kern-0.28ex]}

\opn\dirlim{\underrightarrow{\lim}}
\opn\inivlim{\underleftarrow{\lim}}
\let\to=\rightarrow

\def\Implies{\ifmmode\Longrightarrow \else
     \unskip${}\Longrightarrow{}$\ignorespaces\fi}
\def\implies{\ifmmode\Rightarrow \else
     \unskip${}\Rightarrow{}$\ignorespaces\fi}
\def\iff{\ifmmode\Longleftrightarrow \else
     \unskip${}\Longleftrightarrow{}$\ignorespaces\fi}

\let\:=\colon
\newtheorem{Theorem}{Theorem}[section]
\newtheorem{Lemma}[Theorem]{Lemma}
\newtheorem{Corollary}[Theorem]{Corollary}
\newtheorem{Proposition}[Theorem]{Proposition}
\newtheorem{Remark}[Theorem]{Remark}

\newtheorem{Example}[Theorem]{Example}

\let\epsilon\varepsilon
\let\phi=\varphi
\let\kappa=\varkappa
\def\qed{\ifhmode\textqed\fi
   \ifmmode\ifinner\quad\qedsymbol\else\dispqed\fi\fi}
\def\textqed{\unskip\nobreak\penalty50
    \hskip2em\hbox{}\nobreak\hfil\qedsymbol
    \parfillskip=0pt \finalhyphendemerits=0}
\def\dispqed{\rlap{\qquad\qedsymbol}}

\def\AA{{\cal A}}

\opn\initial{in}
\opn\inim{inm}
\opn\rev{rev}
\opn\Gin{Gin}
\opn\Lex{Lex}
\opn\Shift{Shift}
\opn\shift{shift}
\opn\rate{rate}
\opn\Mon{Mon}
\opn\lex{lex}
\opn\rev{rev}
\opn\red{red}
\opn\max{max}
\opn\min{min}
\opn\initial{in}
\opn\Ker{Ker}
\opn\GL{GL}
\opn\proj{proj}
\opn\Mon{Mon}
\begin{document}

\title{Componentwise linear ideals with minimal or maximal Betti numbers}
\author{J\"urgen Herzog, Takayuki Hibi, Satoshi Murai and Yukihide Takayama }
\address{J\"urgen Herzog, Fachbereich Mathematik und
Informatik, Universit\"at Duisburg-Essen, Campus Essen, 45117
Essen, Germany} \email{juergen.herzog@uni-essen.de}

\address{Takayuki Hibi,
Department of Pure and Applied Mathematics, Graduate School of
Information Science and Technology, Osaka University, Toyonaka,
Osaka 560-0043, Japan} \email{hibi@math.sci.osaka-u.ac.jp}

\address{ Satoshi Murai,
Department of Pure and Applied Mathematics, Graduate School of
Information Science and Technology, Osaka University, Toyonaka,
Osaka 560-0043, Japan} \email{s-murai@ist.osaka-u.ac.jp}

\address{Yukihide Takayama, Department of Mathematical Sciences,
Ritsumeikan University, 1-1-1 Nojihigashi, Kusatsu, Shiga 525-8577,
Japan} \email{takayama@se.ritsumei.ac.jp}
\date{}
\maketitle
\begin{abstract}
We characterize componentwise linear monomial ideals with minimal
Taylor resolution and consider the lower bound for the Betti
numbers of componentwise linear ideals.
\end{abstract}

\bigskip
\bigskip
\bigskip
\bigskip
\bigskip
\section*{Introduction}
 Let $S = K[x_1, \ldots, x_n]$ denote the polynomial ring in $n$
variables over a field $K$ with each $\deg x_i = 1$. Let $I$ be a
monomial ideal of $S$ and $G(I) = \{ u_1, \ldots, u_s \}$ its
unique minimal system of monomial generators. The Taylor
resolution \cite[p.~439]{Eisenbud} provides a graded free
resolution of $S/I$. Fr\"oberg \cite[Proposition 1]{Froeberg}
characterizes the monomial ideals for which the Taylor resolution
is minimal. In most cases it is indeed not minimal, however it
yields the following upper bound for the Betti numbers of $I$.
\[
\beta_i(I)\leq {s\choose i+1}\quad \textnormal{for} \quad
i=0,\ldots,s-1.
\]
This upper bound is reached exactly when the Taylor resolution is
minimal.

On the other hand, Brun and R\"omer \cite[Corollary 4.1]{BR} have
shown that
\[
\beta_i(I)\geq {p\choose i+1}\quad \textnormal{for} \quad
i=0,\ldots,p-1,
\]
where $p$  denotes the projective dimension of $S/I$.

In this note we consider componentwise linear ideals of $S$, and
study the cases when one of the bounds for the Betti numbers
described above is obtained. As one of the main results (Theorem
\ref{Gotzmann}) we have that a componentwise linear monomial ideal
has a minimal Taylor resolution if and only if $|G(I)|\leq n$. We
also prove in Theorem~\ref{generalization} the following result:
assume that $K$ is infinite, and let  $I$ be a componentwise
linear ideal of $S$ and $\Gin(I)$ its generic initial ideal  with
respect to the reverse lexicographic order. Suppose that
$|G(\Gin(I))| = s$, and that $\beta_i(\Gin(I)) = {s \choose i +
1}$ for some $1 \leq i < s$. Then $I$ is Gotzmann and $\beta_i(I)
= {s \choose i + 1}$ for all $i$.

In general, if $I$ is a monomial ideal generated by $s$ elements
and for some $i$ with $0<i<s-1$, the $i$th Betti number of $I$
reaches the Taylor bound, i.e.\,  $\beta_i(I)={s\choose i+1}$,
then, contrary to the previous result,  this does not necessarily
imply that the whole Taylor resolution is minimal, as we show by
examples. While if $\beta_{s-1}(I)\neq 0$, the Taylor resolution
is indeed minimal, as noted by Fr\"oberg \cite[Proposition 1
]{Froeberg}.

Concerning the lower bound we have the following result (Theorem
\ref{lower}):  Let $I\subset \mm^2$ be a componentwise linear
ideal with $\grade I=g$ and $\projdim S/I=p$. Suppose that
$\beta_i(I)={p\choose i+1}$ for some $i$. Then (a) $i\geq g$, and
(b) $\beta_j(I)={p\choose j+1}$ for all $j\geq i$.

\section{Componentwise linear ideals with minimal Taylor
resolution}

A monomial ideal $I$ of $S$ is {\em lexsegment} if, for a monomial
$u$ of $S$ belonging  to $I$ and for a monomial $v$ of $A$ with
$\deg u = \deg v$ and with $v <_{\lex} u$, one has $v \in I$,
where $<_{\lex}$ is the lexicographic order on $S$ induced by the
ordering $x_1 > \cdots > x_n$ of the variables. A lexsegment ideal
$I$ of $S$ is called {\em universal lexsegment} \cite{BNT} if, for
$m = 1, 2, \ldots$, the monomial ideal of the polynomial ring
$K[x_1, \ldots, x_n, x_{n+1}, \ldots, x_{n+m}]$ generated by the
monomials belonging to $G(I)$ is lexsegment. In other words, a
universal lexsegment ideal of $S$ is a lexsegment ideal $I$ of $S$
which remains being lexsegment in the polynomial ring $K[x_1,
\ldots, x_n, x_{n+1}, \ldots, x_{n+m}]$ for $m = 1, 2, \ldots$. It
is known \cite[Corollary 1.4]{MuraiHibi} that a lexsegment ideal
$I$ of $S$ is universal lexsegment if and only if $|G(I)| \leq n$,
where $|G(I)|$ is the number of monomials belonging to $G(I)$.

\begin{Lemma}
\label{universal}
The Taylor resolution of a universal lexsegment ideal
is minimal.
\end{Lemma}

\begin{proof}
Let $I$ be a universal lexsegment ideal of $S$ and $G(I) = \{ u_1,
u_2, \ldots u_s \}$ with $s \leq n$, where $\deg u_1 \leq \deg u_2
\leq \cdots \leq \deg u_s$ and where $u_{i+1} <_{\lex} u_i$ if
$\deg u_i = \deg u_{i+1}$. Then \cite[Lemma 1.2]{MuraiHibi} says
that $m(u_i) = i$, where for a monomial $u$ of $S$ the notation
$m(u)$ stands for the biggest integer $j$ for which $x_j$ divides
$u$. By using Eliahou and Kervaire \cite{E--K} it follows that the
$q$th total Betti number $\beta_q(I)$ of $I$ is $\sum_{i=1}^s {i -
1 \choose q} = {s \choose q + 1}$. Thus $\beta_q(I)$ coincides
with the rank of the $q$th free module of the Taylor resolution of
$I$.
\end{proof}

We refer the reader to \cite{HerzogHibi} for fundamental material
of componentwise linear ideals and Gotzmann ideals. A homogeneous
ideal $I$ of $S$ is {\em Gotzmann} if $\beta_{ij}(I) =
\beta_{ij}(I^{\lex})$ for all $i$ and $j$, where $I^{\lex}$ is the
unique lexsegment ideal of $S$ with the same Hilbert function as
$I$. A homogeneous ideal $I$ of $S$ is {\em componentwise linear}
if $\beta_{ij}(I) = \beta_{ij}(\Gin(I))$ for all $i$ and $j$,
where $\Gin(I)$ is the generic initial ideal of $I$ with respect
to the reverse lexicographic order induced by the ordering $x_1 >
\cdots x_n$ of the variables; see \cite[Theorem 1.1]{AHH} for the
proof of this statement  when $\chara K=0$, and  \cite[Lemma
3.3]{HZ} in the case of a field of arbitrary characteristic.
Examples of componentwise linear ideals are the Gotzmann ideals.

We shall need the following three results  for the proof of the
next theorems.
\begin{Lemma}
\label{easy} A componentwise  ideal $I$  is Gotzmann if and only
if $\Gin(I)$ is Gotzmann.
\end{Lemma}

\begin{proof} Since $I$ is componentwise linear we have
$\beta_{ij}(I)=\beta_{ij}(\Gin(I)$ for all $i$ and $j$. On the
other hand, $I^{\lex}=\Gin(I)^{\lex}$. Therefore,
$\beta_{ij}(I)=\beta_{ij}(I^{\lex})$ if and only if
$\beta_{ij}(\Gin(I))=\beta_{ij}(\Gin(I)^{\lex})$. This proves the
assertion.
\end{proof}

Next we have

\begin{Lemma}
\label{nice} Let $I$ be a stable ideal and $u \in G(I)$ with $m(u)
= i \geq 2$.  Then there exists  $w \in G(I)$ with $\deg w \leq
\deg u$ and with $m(w) = i - 1$. In particular, $\max\{m(u)\:\;
u\in G(I)\}\leq |G(I)|$.
\end{Lemma}

\begin{proof}
Let $u = v x_i^N$ with $m(v) \leq i - 1$ and with $N \geq 1$.
Since $I$ is stable, one has $v x_{i-1}^N \in I$. Thus there is $w
\in G(I)$ which divides $v x_{i-1}^N$. Since $u \in G(I)$, it
follows that $w$ cannot divide $v$. Hence $m(w) = i - 1$ and $\deg
w \leq \deg u$.
\end{proof}

\begin{Lemma}
\label{HibiMurai} Let $I$ be a stable ideal with $\max\{m(u)\:\;
u\in G(I)\}=|G(I)|$. Then $I$ is a Gotzmann ideal with $|G(I)|\leq
n$. \end{Lemma}

\begin{proof}By using Lemma \ref{nice}
we may assume that $G(I) = \{ u_1, \ldots, u_s \}$, where $\deg
u_1 \leq \cdots \leq \deg u_s$ and where $m(u_i) = i$ for $1 \leq
i \leq s$. Now, \cite[Lemma 1.4]{MuraiHibi} guarantees that there
exists a universal lexsegment ideal $L$ of $S$ with $G(L) = \{
w_1, \ldots, w_s \}$ such that $\deg u_i = \deg w_i$ for $1 \leq i
\leq s$. Again, by \cite[Lemma 1.2]{MuraiHibi}, one has $m(w_i) =
i$ for $1 \leq i \leq s$. Hence the Eliahou--Kervaire formula
\cite{E--K} implies  that $\beta_{ij}(I) = \beta_{ij}(L)$ for all
$i$ and $j$. In particular $I$ and $L$ have the same Hilbert
function. Thus $L = I^{\lex}$. It follows that $I$ is Gotzmann. Of
course, $|G(I)| \leq n$, since $m(u)\leq n$ for all $u\in G(I)$.
\end{proof}

The equivalence of the statements (a) and (b) in the following
Theorem  has been proved for the special case of stable ideals in
the paper \cite{OT}.

\begin{Theorem}
\label{Gotzmann} Let $I$ be a componentwise linear monomial ideal
of $S$. Then the following conditions are equivalent:
\begin{enumerate}
\item[(a)] The Taylor resolution of $I$ is minimal.

\item[(b)] $\max\{m(u)\:\; u\in G(I)\}=|G(I)|$.

\item[(c)]  $I$ is Gotzmann with $|G(I)| \leq n$.

\end{enumerate}
\end{Theorem}

\begin{proof} (a) \implies (b):  Let $|G(I)| = s$ and  let $J = \Gin(I)$. Then $J$ is
strongly stable with $|G(J)| = s$, and the Taylor resolution of
$J$ is minimal as well since the Betti numbers do not change.
Since $\beta_{s-1}(J) \neq 0$, the Eliahou--Kervaire formula
implies that there exists   a monomial $u \in G(J)$ with $m(u) =
s$.

(b) \implies (c):  By using  Lemma \ref{HibiMurai}, it follows
that $J$ is a Gotzmann ideal with $|G(J)|\leq n$.  Since $J^{\lex}
= I^{\lex}$ and $J$ is Gotzmann, we have
 $\beta_{ij}(I) = \beta_{ij}(J) = \beta_{ij}(J^{\lex})= \beta_{ij}(I^{\lex})$ for all
$i$ and $j$. Thus  $I$ is Gotzmann with $|G(I)| \leq n$.

(c)\implies (a): Since $|G(I^{\lex})| = |G(I)|\leq n$, the
lexsegment ideal $I^{\lex}$ is a universal lexsegment ideal. Since
$\beta_i(I) = \beta_i(I^{\lex})$ for all $i$, and since by Lemma
\ref{universal} the Taylor resolution of $I^{\lex}$ is minimal, it
follows that the Taylor resolution of $I$ is minimal, as required.
\end{proof}

\begin{Remark}{\em  If some of the Betti numbers of a monomial ideal
reaches the Taylor bound, then this does not necessarily imply
that the whole Taylor resolution is minimal. Given integers $2\leq
i<s-2$, we consider the ideal $I$ generated by $x_1y_1,x_2y_2,
\ldots,x_{s-1}y_{s-1}, y_1\cdots y_{i}$. Then it is easily checked
that $\beta_j(I)={s\choose j+1}$ for $j=0,\ldots i-2$ and
$\beta_j(I)<{s\choose j+1}$ for $i-2<j\leq s-1$.

On the other hand, if $I$ is a monomial ideal with
$G(I)=\{u_1,u_2,\ldots,u_s\}$, and for some $0\leq i\leq s$ we
have $\beta_i(I)={s\choose i+1}$, then $\beta_j(I)={s\choose j+1}$
for all $j\leq i$. Indeed, if $\beta_i(I) = {s\choose i+1}$. Then
from the definition of the Taylor complex it follows that for an
arbitrary subset $\{u_{j_1},u_{j_2},\ldots, u_{j_i}\}$ of $G(I)$
of cardinality $i$, one has that $\lcm(u_{j_1},\dots, ,u_{j_i})
\neq \lcm(u_{j_1},\ldots, u_{j_{k-1}},u_{j_{k+1}},\ldots,u_{j_i})$
for all $k=1,...,i$. It obvious that similar inequalities hold for
any subset of $\{u_{j_1},u_{j_2},\ldots, u_{j_i}\}$. This clearly
implies that $\beta_j(I) = {s\choose j+1}$  for all $j\leq i$.}
\end{Remark}

\begin{Theorem}
\label{generalization} Let $I$ be a componentwise linear ideal of
$S$ with $|G(\Gin(I))| = s$.  Suppose  that $\beta_i(\Gin(I)) = {s
\choose i + 1}$ for some $1 \leq i < s$. Then $I$ is Gotzmann and
$\beta_i(I) = {s \choose i + 1}$ for all $i$.
\end{Theorem}

\begin{proof}
Let $q \geq 1$
be the biggest integer for which
$m(u) = q$
for some $u \in G(\Gin(I))$.
Lemma \ref{nice} says that
for each $j \leq q$ there is
$u \in G(\Gin(I))$ with
$m(u) = j$.
Fix $i_0$ with
$\beta_{i_0}(\Gin(I)) = {s \choose i_0 + 1}$.
Since
\begin{eqnarray*}
\beta_{i_0}(\Gin(I))
& = & \sum_{u \in G(\Gin(I))} {m(u) - 1 \choose i_0}
\\
& = & \sum_{u \in G(\Gin(I)), \, m(u) > i_0}
{m(u) - 1 \choose i_0}
\\
& \leq &
{i_0 \choose i_0} + {i_0 + 1 \choose i_0} +
\cdots + { q - 2 \choose i_0}
+ (s - q + 1){q - 1 \choose i_0}
\\
& \leq &
{i_0 \choose i_0} + {i_0 + 1 \choose i_0} +
\cdots + { s - 1 \choose i_0}
= {s \choose i_0 + 1},
\end{eqnarray*}
it follows that $q = s$. Hence Lemma \ref{HibiMurai} implies that
$\Gin(I)$ is Gotzmann with $|G(\Gin(I))|\leq n$. Therefore Theorem
\ref{Gotzmann} guarantees that $\beta_i(\Gin(I)) = {s \choose i +
1}$ for all $i$. Since $I$ is componentwise linear, we have
$\beta_i(I)=\beta_i(\Gin(I))$ for all $i$, and by Lemma \ref{easy}
that $I$ itself is Gotzmann, as desired.
\end{proof}

\section{Componentwise linear ideals with minimal Betti numbers}

It is not surprising that a stable  monomial ideal is rarely a
complete intersection. One such example is
$I=(x_1,\ldots,x_{s-1},x_s^d)$ where $1\leq s\leq n$ and $d\geq
1$. If $I$ is contained in the square of the graded maximal ideal
$\mm=(x_1,\ldots,x_n)$ of $S$, as we may always assume, then there
exist even less such ideals. Indeed we have

\begin{Lemma}
\label{less} Let  $I\subset \mm^2$ be a stable ideal. If $I$ is a
complete intersection, then $I=(x_1^d)$ for some $d\geq 2$.
\end{Lemma}

\begin{proof} Let $G(I)=\{u_1,\ldots, u_s\}$. We may assume that
$m(u_s)\geq m(u_{s-1})\geq \cdots \geq m(u_1)$. Assume $s\geq 2$.
Then $(u_1,u_2)$ is stable and a complete intersection. It is
clear that $u_1=x_1^{d}$ for some $d\geq 2$. By  Lemma \ref{nice}
we have  $m(u_{2})=2$ and $\deg u_2\geq d$. Since $u_1$ and $u_2$
are coprime, it follows that $u_2=x_2^c$ with $c\geq d$. Since $I$
is stable, $x_1x_2^{c-1}\in (u_1,u_2)$, a contradiction.
\end{proof}

Now we can show

\begin{Theorem}
\label{lower} Let $I\subset \mm^2$ be a componentwise linear ideal
with $\grade I=g$ and\\ $\projdim S/I=p$. Suppose that
$\beta_i(I)={p\choose i+1}$ for some $i$. Then
\begin{enumerate}
\item[(a)] $i\geq g$.

\item[(b)] $\beta_j(I)={p\choose j+1}$ for all $j\geq i$.
\end{enumerate}
\end{Theorem}

\begin{proof} We may replace $I$ by $\Gin(I)$ and hence may assume that $I$ is a stable ideal.
As usual we set $m_i(I)=|\{u\in G(I)\:\; m(u)=i\}|$. We first show
that $m_i(I)>1$ for $i=2,\ldots, g$. Indeed, since $I$ is a stable
monomial ideal of grade $g$, the sequence
$x_n,x_{n-1},\ldots,x_{n-g+1}$ is a system of homogeneous
parameters for the standard graded $K$-algebra $S/I$. Hence $S/I$
modulo the sequence $x_n,x_{n-1},\ldots,x_{n-g+1}$ is of dimension
$0$ and isomorphic to $K[x_1,\ldots, x_g]/J$ where $J$ is again a
stable ideal with $m_i(J)=m_i(I)$ for $i=1,\ldots,g$. We may
assume that $g>1$. Then $m_g(I)=m_g(J)=\dim_K(J\: x_g)/J$. Since
$J$ is stable, $J\: x_g=J\: (x_1,\ldots,x_g)$, and it follows that
$m_i(I)=\dim_KJ\: (x_1,\ldots,x_g)/J>1$, since  $S/J$ is not
Gorenstein. In fact, if $S/J$ would be Gorenstein, then, since
$\dim S/J=0$ and $J$ is a monomial ideal, it would follow that $J$
is a complete intersection, contradicting  Lemma \ref{less}. If
$g>2$, then we consider $K[x_1,\ldots, x_g]/J$ modulo $x_g$ and
repeat the argument. This can be done as long as $g>2$.

Now we use the Eliahou--Kervaire formula and get
\begin{eqnarray*}
\beta_{i}(I) & = & \sum_{u \in G(I)} {m(u) - 1 \choose
i}=\sum_{j=i}^{p-1}m_{j+1}(I){j\choose i}\geq
\sum_{j=i}^{p-1}{j\choose i}={p\choose i+1}
\end{eqnarray*}
with equality if and only if all $m_j(I)=1$ for $j=i+1,\ldots,p$.

For the proof of  (a) we may assume that $g\geq 2$. Then
$m_g(I)\geq 2$, and we can have the equality
$\beta_{i}(I)={p\choose i+1}$ only if $i+1> g$.

On the other hand, if $\beta_{i}(I)={p\choose i+1}$ for some
$i\geq g$. Then $m_{j}(I)=1$ for all $j=i+1,\ldots,p$, and hence
the Eliahou--Kervaire formula implies that $\beta_{j}(I)={p\choose
j+1}$ for all $j\geq i$, as desired.
\end{proof}

In view of the preceding theorem  and in view of Lemma
\ref{universal} the following consequence is immediate.

\begin{Corollary}
\label{g} Let $I\subset \mm^2$ be a componentwise linear ideal
whose Taylor resolution is minimal. Then $\grade I=1$. In
particular, any universal lexsegment ideal $I\subset \mm^2$  is of
the form $x_1J$ where $J$ is a monomial ideal.
\end{Corollary}


\begin{thebibliography}{99}

\bibitem{AHH} A. Aramova, J. Herzog, T.  Hibi,  Ideals with
stable Betti numbers, {\em Adv. Math.}  {\bf 152} no. 1  (2000),
72--77.



\bibitem{BNT}
E. Babson, I. Novik and R. R. Thomas,
Reverse lexicographic and lexicographic shifting,
{\em J. Algebraic Combin.}, to appear.


\bibitem{BR} M. Brun and T. R\"omer,  Betti numbers of ${\Bbb Z}\sp
n$-graded modules,  {\em Comm. Algebra} {\bf 32} no. 12  (2004),
4589--4599.


\bibitem{Eisenbud}
D. Eisenbud,
``Commutative Algebra with a View Toward
Algebraic Geometry,'' Springer--Verlag, 1995.

\bibitem{E--K}
S. Eliahou and M. Kervaire, Minimal resolutions of some monomial
ideals, {\em J. of Algebra} {\bf 129} (1990), 1 -- 25.

\bibitem{Froeberg} R.\ Fr\"oberg, Some complex constructions with applications to Poincaré series.
Séminaire d'Algèbre Paul Dubreil 31\`{e}me ann\'{e} (Paris,
1977--1978), pp.\ 272--284, Lecture Notes in Math., {\bf 740},
Springer, Berlin, 1979.


\bibitem{HerzogHibi}
J. Herzog and T. Hibi,
Componentwise linear ideals,
{\em Nagoya Math. J.} {\bf 153} (1999),
141 -- 153.

\bibitem{HZ} J. Herzog, X.  Zheng, Notes on the multiplicity
conjecture, {\em Collect. Math. } {\bf 57}  no. 2  (2006),
211--226.




\bibitem{MuraiHibi}
S. Murai and T. Hibi,
The depth of an ideal with a given Hilbert function,
preprint, 2006,
arXiv:math.AC/0608188.


\bibitem{OT}  M. Okudaira and Y. Takayama, Monomial ideals with linear
quotients whose Taylor resolutions are minimal, preprint, 2006,
arXiv:math.AC/0610168.



\end{thebibliography}
\end{document}